\documentclass[11pt]{article}

\usepackage[centertags]{amsmath}
\usepackage{amsfonts}
\usepackage{amssymb}
\usepackage{theorem}
\usepackage{epsfig}
\usepackage{color}
\theoremstyle{break} \newtheorem{theorem}{Theorem}[section]
\theoremstyle{break} 
\theoremstyle{break}        
\theoremstyle{break} \newtheorem{lemma}[theorem]{Lemma}
\theoremstyle{break} \newtheorem{corollary}[theorem]{Corollary}
\theoremstyle{break} 
\theoremstyle{break} 
\theoremstyle{break}
{\theorembodyfont{\rmfamily}\newtheorem{remark}[theorem]{Remark}}
{\theorembodyfont{\rmfamily}}
\theoremstyle{break} 
\theoremstyle{break} 
\theoremstyle{break} 
\theoremstyle{break} 
\numberwithin{equation}{section}

\newcommand{\R}{{\mathbb{R}}}
\newcommand{\D}{{\mathbb{D}}}
\newcommand{\C}{{\mathbb{C}}}

\renewcommand{\P}{{\mathbb{P}}}
\renewcommand{\H}{{\mathbb{H}}}

\def\Re{\mathop{{\rm Re}}}
\def\Im{\mathop{{\rm Im}}}

\def\arcsinh{\mathop{{\rm arcsinh}}}

\begin{document}

\renewcommand{\thefootnote}{}
\stepcounter{footnote}
\begin{center}
{\bf \Large Metrics with conical singularities on the  sphere
 and sharp\\[2mm] extensions of the theorems of Landau and Schottky
\footnote{2000 Mathematics Subject
 Classification: Primary  30F45}}
\end{center}
\renewcommand{\thefootnote}{\arabic{footnote}}
\setcounter{footnote}{0}
\begin{center}
{\large Daniela Kraus, Oliver Roth and Toshiyuki Sugawa}\\[2mm]

\end{center}


\bigskip
\renewcommand{\thefootnote}{\arabic{footnote}}
\smallskip
\begin{center}
\begin{minipage}{13cm} {\bf Abstract.}   
An explicit formula for the generalized hyperbolic metric on the thrice--punctured
sphere $\P \backslash \{z_1, z_2, z_3\}$ with singularities of order $\alpha_j \le 1$ at $z_j$ is obtained in all possible cases
$\alpha_1+\alpha_2+\alpha_3 >2$. The existence and uniqueness of such a metric
was proved long time ago by Picard \cite{Pic1905} and Heins \cite{Hei62}, while
explicit formulas for the cases $\alpha_1=\alpha_2=1$  were given earlier by
Agard \cite{AG} and recently by Anderson, Sugawa, Vamanamurthy and Vuorinen
\cite{A}. We also establish precise and explicit lower bounds for the generalized
hyperbolic metric. This extends work of Hempel \cite{Hem79} and Minda
\cite{Min87b}.  As applications, sharp versions of Landau-- and Schottky--type
theorems for meromorphic functions are obtained.
\end{minipage}
\end{center}
\footnotetext{D.K.~and O.R.~were supported
by a DFG grant (RO 3462/3--1). 
 T.S.~was supported in part by JSPS Grant-in-Aid for Scientific
Research (B), 17340039 and for Exploratory Research, 19654027.
\hfill{To~appear:~{\it Math.~Z.}}}
\section{Introduction}

Let $\P$ denote the Riemann sphere 
endowed with its canonical complex structure and let $\Omega \subseteq \P$
be a subdomain.
We say a  conformal Riemannian metric $\lambda(z) \, |dz|$ on 
$\Omega \backslash \{ p \}$ 
has a singularity of order $ \alpha \le 1$ at the
point $p \in \Omega$, if, in local coordinates,
$$ \log \lambda(z)  = \begin{cases}
 -\alpha \log |z-p|+O(1) & \text{ if } \alpha <1 \, \\ \\
- \log |z-p|-\log \left( -\log |z-p| \right)+O(1) & \text{ if } \alpha=1
\, \end{cases}$$
 as $z \to p$. 
Geometrically, the singular surface $(\Omega,\lambda(z) \, |dz|)$ 
looks like an ice--cream cone at $p$ if  $\alpha<1$.
 If $\alpha=1$, then $(\Omega,\lambda(z)\, |dz|)$ has a cusp at $p$.
We therefore call $p$ a {\it conical singularity} or {\it corner of order} $\alpha$
if  $\alpha<1$ and a {\it cusp}  if  $\alpha=1$.
It is also customary to say that a conformal Riemannian metric with a 
conical singularity $p$ of order 
$\alpha<1$ has the angle $2\pi(1-\alpha)$ at the point $p$, see \cite{A}.
Conical singularities are  typical for  conformal metrics. For instance, if the curvature of $\lambda(z) \, |dz|$ is bounded below and above by
 negative constants, then $\lambda(z) \, |dz|$ only has corners 
or cusps as point
 singularities  (see \cite{Hei62,McO93,KR2007}). For nonnegatively curved metrics
with finite energy  only corners occur (see \cite{Yun03,KR2008}).

\medskip

It is well--known (see \cite{Pic1905,Hei62}) that for
 $n \ge 3$ distinct points $ z_1 , \ldots \, , z_n \in \P$ and 
real parameters $ \alpha_1, \ldots,
\alpha_n  \in (-\infty,1]$ 
there exists a  conformal
 metric on the $n$--punctured sphere $\P\backslash \{z_1, \ldots, z_n\}$ with constant curvature $-1$ and singularities of order $\alpha_j$ 
at $z_j$ if and only if
\begin{equation} \label{eq:gaussbonnet0}
\sum \limits_{j=1}^n
\alpha_j >2 \, .
\end{equation}
In this case, this metric is uniquely determined and will be called 
{\it generalized hyperbolic metric}
with singularities of order $\alpha_j$ at $z_j$. 

\smallskip

We note that the necessity of  condition (\ref{eq:gaussbonnet0}) comes from the Gauss--Bonnet theorem; the sufficiency is a special case of the classical Schwarz--Picard problem which has been solved by Picard \cite{Pic1905} and Heins \cite{Hei62}, see also 
Bieberbach \cite{Bie16}, McOwen \cite{McO88,McO93} and Troyanov
\cite{Tro90}.  
The terminology {\it generalized hyperbolic metric} is motivated by the fact that  
 if  all singularities are cusps, then one gets back
the standard hyperbolic metric on the punctured sphere $\P \backslash \{ z_1, \ldots, z_n\}$.

\smallskip


\medskip

We are primarily interested in the case of the thrice--punctured sphere
$\P \backslash \{ z_1, z_2, z_3\}$.  Note that in this case
\begin{equation} \label{eq:gaussbonnet}
 0 < \alpha_1 \le 1 \, , \quad 0 < \alpha_2 \le 1 \, , \quad 0 <
\alpha_3 \le 1 \, , \quad \alpha_1+\alpha_2+\alpha_3>2 \, .
\end{equation}
Using a  M\"obius transformation, which sends $z_1$ to $0$, $z_2$ to $1$
and $z_3$ to $\infty$, 
we may henceforth assume  that $z_1=0$, $z_2=1$ and $z_3=\infty$ and
shall  denote by
$\lambda_{\alpha_1,\alpha_2,\alpha_3}(z)\, |dz|
$
the generalized hyperbolic metric with conical singularities of order
$\alpha_1$, $\alpha_2$ and $\alpha_3$ at $z_1=0$, $z_2=1$ and $z_3= \infty$.
In this situation, the Riemannian metric
$\lambda_{\alpha_1,\alpha_2,\alpha_3}(z) \, |dz|$ on $\P \backslash
\{0,1,\infty\}$ can be
described in terms of a single density 
function $\lambda_{\alpha_1,\alpha_2,\alpha_3}$ defined on the twice--punctured
plane $\C'':=\C \backslash \{
0,1\}$ (see Section \ref{sec:preliminaries} below for details). We call $\lambda_{\alpha_1,\alpha_2,\alpha_3}$
the {\it generalized hyperbolic density} of order
$(\alpha_1,\alpha_2,\alpha_3)$ on $\C''$.

\medskip

Explicit and very useful formulas for
$\lambda_{1,1,\alpha_3}(z) $ have been obtained by Agard \cite{AG}
for $\alpha_3=1$ and recently by
Anderson, Sugawa, Vamanamurthy and Vuorinen
\cite{A} for $\alpha_3 \in (0,1]$. Hempel \cite{Hem79} (see also Minda
\cite{Min87b}) proved a sharp, explicit and easy--to--use lower bound for the
standard hyperbolic density $\lambda_{1,1,1}(z)$. In combination with
Agard's  formula for $\lambda_{1,1,1}(z)$ this has led to
precise bounds in the classical theorems of Landau and Schottky\footnote{We refer the reader to the monographs \cite{Burckel} and \cite{Hayman} for an
introduction to Landau's and Schottky's theorem and to the recent paper \cite{AQV}
for connections  of Schottky's theorem
with quasiconformal maps and  modular equations.}
for analytic functions in the open unit disk $\D:=\{z \in \C \, : \, |z|<1\}$ omitting the values $0$ and $1$
(see Ahlfors \cite{Ahl38}, Hayman \cite{Hay},
Hempel \cite{Hem79,Hem80}, Jenkins \cite{Jen} and e.g.~Li \& Qi \cite{LQ2007}).

\medskip

In this note, we extend the above results
to the {\it generalized} 
hyperbolic metric and provide sharp extensions of theorems of Landau and 
Schottky type for meromorphic functions not necessarily omitting the
values $0$, $1$ and $\infty$.

\medskip
 
The  paper is organized in the following way.
The main results are described and discussed in Section \ref{sec:results}.
Section \ref{sec:preliminaries}  contains  a
 quick review of the necessary background material about
conformal metrics, while  Section \ref{sec:proofs} is devoted to the proofs of
the results.
%
%
We start in \S \ref{sec:formulas} with Theorem \ref{thm}, which provides
an explicit formula for the generalized hyperbolic density
$\lambda_{\alpha_1,\alpha_2,\alpha_3}(z) $
in all possible cases (see (\ref{eq:gaussbonnet})). This  generalizes the
results of Agard \cite{AG} and Anderson, Sugawa, Vamanamurthy and Vuorinen
\cite{A}, which are easily seen to be special cases of Theorem \ref{thm}, to
the most general situation. Our method of proof differs from that in 
\cite{AG,A} as we base our proof
on Liouville's representation formula (Theorem \ref{thm:liouville}) for
constantly curved conformal Riemannian metrics. The use of Liouville's theorem
will also facilitate proving sharpness of (most of) our results.
 In \S \ref{sec:mon} we give
a {\it sharp} and  {\it explicit} lower bound for the generalized hyperbolic
metric, see Theorem \ref{thm:lowerbound}. This extends
 the earlier work of
Hempel \cite{Hem79} and Minda \cite{Min87b}, which deals with the special case
of the standard hyperbolic metric, to the generalized hyperbolic metric.
Our method is based on a new device, which we call
the Gluing lemma (see Lemma \ref{glueing}) and which allows a rather 
quick proof of Theorem \ref{thm:lowerbound}. 

\medskip

These new information about the generalized hyperbolic metric, which are
perhaps also
interesting in their own right,
  are then applied to study value distribution properties of functions
  meromorphic in the unit disk. For that purpose it is sufficient to consider
 the cases $\alpha_1=1-1/j$, $\alpha_2=1-1/k$,
$\alpha_3=1-1/l$,  where
$j,k,l \ge 2$ are integers (or $=\infty)$ such that according to 
(\ref{eq:gaussbonnet})
$$ \frac{1}{j}+\frac{1}{k}+\frac{1}{l}<1 \, . $$
In this way, we are led to {\it sharp} extensions of the theorems of Landau and
Schottky for meromorphic functions belonging to the classes
\begin{equation} \label{eq:eq}
\hspace*{-0.4cm}\begin{array}{rl}
{\cal M}_{j,k,l}& :=\{f \text{ meromorphic in } \D \text{ such that (i) all
  zeros of } f \text{ have
 order } \ge j, \text{ (ii) all}\\ & \hspace{1cm}\text{ zeros}   \text{ of } f-1  \text{ have order } \ge k
\text{ and }  \text{(iii) all
poles of } f \text{ have order }\ge l\}.
\end{array}
\end{equation}
These results, which  are discussed in Paragraph \ref{sec:apps} and proved
in Section \ref{sec:proofs},
generalize the results in \cite{Ahl38,Hay,Hem79,Hem80,Jen,LQ2007}, which deal
with the particular case of analytic functions in $\D$ omitting the values $0$
and $1$, i.e., the class ${\cal M}_{\infty,\infty,\infty}$, to the much wider
classes ${\cal M}_{j,k,l}$.

\section{Results} \label{sec:results}

\subsection{Explicit formulas} \label{sec:formulas}

The explicit formula for the generalized hyperbolic density
$\lambda_{\alpha_1,\alpha_2,\alpha_3}(z)$, 
 which will be stated momentarily, is necessarily a bit
technical, so we first need to introduce some notation.
Let $\alpha_1, \alpha_2, \alpha_3$ be real parameters satisfying
condition (\ref{eq:gaussbonnet}). We define
\begin{equation} \label{eq:hyp1}
\alpha := \displaystyle \frac{\alpha_1+\alpha_2-\alpha_3}{2} \, , \qquad
\beta := \displaystyle \frac{\alpha_1+\alpha_2+\alpha_3-2}{2} \, , \qquad 
\gamma := \alpha_1 \, . 
\end{equation}
Then $0 < \beta \le \alpha$ and $\alpha+\beta \le \gamma \le 1$. We also
consider the hypergeometric functions
$$
 \varphi_1(z):=  \displaystyle F(\alpha,\beta,\gamma;z) \, , \qquad 
\varphi_2(z):=  \displaystyle F(\alpha,\beta,\alpha+\beta-\gamma+1;1-z) \, .
$$
Note that $\varphi_1$ is analytic in $\C \backslash [1,+\infty)$ and
$\varphi_2$ is analytic in $\C \backslash [-\infty,0]$.

\begin{theorem}[Corners at $\mathbf{z=0}$ and $\mathbf{z=1}$] \label{thm}
Let $0< \alpha_1, \alpha_2<1$ and $0 < \alpha_3 \le 1$ such that $\alpha_1+\alpha_2+\alpha_3 >2$. Then
$$\lambda_{\alpha_1,\alpha_2,\alpha_3}(z)=\frac{1}{|z|^{\alpha_1} |1-z|^{\alpha_2}}
\frac{2 \, K_3}{K_1 |\varphi_1(z)|^2+K_2 |\varphi_2(z)|^2+2 \Re (\varphi_1(z) \varphi_2(\overline{z}))} \,
,$$
where
\begin{equation} \label{eq:consts}
\begin{array}{lcl}
K_1&:=& \displaystyle -\frac{\Gamma (\gamma -\alpha ) \, \Gamma (\gamma -\beta )}{\Gamma (\gamma ) \, \Gamma
   (\gamma-\alpha -\beta )} \, , \qquad 
K_2 :=  \displaystyle-\frac{\Gamma (\alpha +1-\gamma )\,  \Gamma (\beta+1 -\gamma )}{\Gamma (1-\gamma ) \,
   \Gamma (\alpha +\beta +1-\gamma )}\, , \\[6mm]
K_3 &:=&  \displaystyle\sqrt{\frac{\sin (\pi \alpha) \, \sin(\pi \beta)}{\sin
    ( \pi(\gamma-\alpha)) \, \sin (\pi (\gamma-\beta))}}\cdot
\frac{\Gamma(\alpha+\beta+1-\gamma) \,  \Gamma(\gamma)}{\Gamma(\alpha) \,
  \Gamma(\beta)} \, , 
   \end{array}
\end{equation}
and $\alpha,\beta, \gamma$ are defined as in (\ref{eq:hyp1}).
\end{theorem}

The previously known formulas for $\lambda_{1,1,1}(z)$ (see
\cite{AG}) and for $\lambda_{1,1,\alpha_3}(z)$, $0 <\alpha_3 \le 1$,
(see \cite{A}) can easily be obtained
from Theorem \ref{thm} by letting $\alpha_j \nearrow 1$ for $j=1,2$.
We omit the details.

\subsection{Sharp lower bounds} \label{sec:mon}

The aim of this section is to provide a sharp lower bound
for $\lambda_{\alpha_1,\alpha_2,\alpha_3}(z)$.
This generalizes the work of Hempel \cite{Hem79} and
Minda \cite{Min87b}, who give a precise lower bound for 
$\lambda_{1,1,1}(z)$,
 to the  general case (\ref{eq:gaussbonnet}).

\begin{theorem}[A sharp explicit lower bound] \label{thm:lowerbound}
Let $\alpha_1,\alpha_2,\alpha_3$ be real parameters satisfying condition
(\ref{eq:gaussbonnet}) and let
\begin{equation} \label{eq:12}
 C_1:=\frac{1}{1-\alpha_1} \arcsinh \left(
  \frac{1-\alpha_1}{\lambda_{\alpha_1,\alpha_2,\alpha_3}(-1)} \right) \, , \qquad
   C_3:=\frac{1}{1-\alpha_3} \arcsinh \left(
  \frac{1-\alpha_3}{\lambda_{\alpha_1,\alpha_2,\alpha_3}(-1)} \right) \, .
\end{equation}
Then
\begin{equation}\label{eq:lowerbound} 
\lambda_{\alpha_1,\alpha_2,\alpha_3}(z) \ge \begin{cases} \displaystyle \frac{1-\alpha_1}{|z| \sinh
    \big[ (1-\alpha_1) \left( C_1-\log|z| \right) \big]}& \hspace{0.8cm}
  |z| \le 1\, ,
  \, z\not=0,1 \, \\
& \text{ if } \\
\displaystyle \frac{1-\alpha_3}{|z| \sinh \big[ (1-\alpha_3)\left( C_3+\log|z| \right) \big]
} & \hspace{0.8cm}|z|>1 \, . 
\end{cases}
\end{equation}
Equality holds if and only if $z=-1$.
\end{theorem}

\begin{remark}[Limit cases $\alpha_1 \nearrow 1$ and $\alpha_2 \nearrow 1$]
If $\alpha_1=1$ and/or $\alpha_3=1$, then the formulas for $C_1$ and
$C_3$ as well as the lower bounds for $\lambda_{\alpha_1,\alpha_2,\alpha_3}$ are to be understood
in the limit sense $\lim_{\alpha_1 \to 1-}$ resp.~$\lim_{\alpha_3 \to 1-}$.
\end{remark}

\begin{remark}[Computation of $\lambda_{\alpha_1,\alpha_2,\alpha_3}(-1)$]
The sharp lower bound (\ref{eq:lowerbound}) for the generalized hyperbolic density requires
the computation of the particular value
$\lambda_{\alpha_1,\alpha_2,\alpha_3}(-1)$.
Using the M\"obius transformation $T(z)=z/(z-1)$, which fixes $z=0$ and
interchanges $z=1$ with $z=\infty$,
 and the easily verified fact
$
\lambda_{\alpha_1,\alpha_2,\alpha_3}(z)=\lambda_{\alpha_1,\alpha_3,\alpha_2}(T(z))
\, |T'(z)|$, 
we get 
$$
\lambda_{\alpha_1,\alpha_2,\alpha_3}(-1)=\frac{\lambda_{\alpha_1,\alpha_3,\alpha_2}(1/2)}{4}
\, .$$
In view of Theorem \ref{thm}, the computation of
$\lambda_{\alpha_1,\alpha_2,\alpha_3}(-1)$ is thereby essentially reduced to the evaluation of 
two hypergeometric functions  at the point $z=1/2$,
which can effectively be achieved using the rapidly converging hypergeometric
series.
\end{remark}

We wish to single out the special case $\alpha_3=\alpha_1$ of Theorem
\ref{thm:lowerbound}, because then 
 not only $C_1=C_3$ holds, but also the value of 
$\lambda_{\alpha_1,\alpha_2,\alpha_1}(-1)$ can explicitly be computed 
 in terms of the Gamma function.

\begin{corollary}[The case $\alpha_3=\alpha_1$] \label{cor:lowerbound}
Let $\alpha_1,\alpha_2 \in (0,1]$ such that $2 \alpha_1+\alpha_2>2$.
 Then
\begin{equation} \label{eq:f}
 \lambda_{\alpha_1,\alpha_2,\alpha_1}(-1)=2 \sqrt{\frac{\tan \left( \frac{\pi}{2} \left(
        \frac{\alpha_2}{2}+\alpha_1 \right) \right)}{\tan \left( \frac{\pi}{2} \left(
        \frac{\alpha_2}{2}-\alpha_1 \right) \right)}} \cdot 
\frac{\Gamma \left( \frac{\alpha_1}{2}-\frac{\alpha_2}{4}+\frac{1}{2}\right)
  \Gamma \left(\frac{\alpha_1}{2}+\frac{\alpha_2}{4} \right)}{\Gamma \left(\frac{\alpha_1}{2}-\frac{\alpha_2}{4}\right) \Gamma \left(\frac{\alpha_1}{2}+\frac{\alpha_2}{4}-\frac{1}{2} \right)}
 \end{equation}
and
$$ \lambda_{\alpha_1,\alpha_2,\alpha_1}(z) \ge \frac{1-\alpha_1}{|z| \sinh \big[ (1-\alpha_1)\left( C_1+\big|\log|z|\big| \right) \big]}$$
for all $z \in \C''$ with equality if and only if $z=-1$.
Here, $C_1$ is given by (\ref{eq:12}) with $\alpha_3=\alpha_1$.
\end{corollary}

For $\alpha_1=\alpha_2=\alpha_3=1$, Corollary \ref{cor:lowerbound}
 further reduces to the sharp bound
$$ \lambda_{1,1,1}(z) \ge \frac{1}{|z| \left( 1/\lambda_{1,1,1}(-1)+\big|\log|z| \big|
  \right)} \quad \text{ with } \quad \lambda_{1,1,1}(-1)=\frac{4
  \pi^2}{\Gamma(1/4)^4} \approx 0.228473 \, ,$$
which was first proved by Hempel \cite{Hem79} (see also Minda
\cite{Min87b}).


\subsection{Applications: Theorems of Landau and Schottky type}
\label{sec:apps}

 Hempel \cite{Hem79,Hem80} (see also Jenkins \cite{Jen}), 
Minda \cite{Min87b}
and Li \& Qi \cite{LQ2007} proved sharp Landau and Schottky type
theorems for functions in  ${\cal
  M}_{\infty,\infty,\infty}$, i.e.,
analytic functions in $\D$ omitting $0$ and~$1$,
  with the help of the standard hyperbolic metric
$\lambda_{1,1,1}(z) \, |dz|$. Using the explicit formula and the sharp lower 
bounds for the generalized hyperbolic metric obtained in the previous
sections we now generalize these results by proving {\it sharp}
versions of Landau and Schottky type theorems
for functions belonging to  the much larger classes ${\cal M}_{j,k,l}$ (see (\ref{eq:eq})).
Here $j,k,l \ge 2$ are integers (or $=+\infty$) such that 
$$ \frac{1}{j}+\frac{1}{k}+\frac{1}{l}<1 \, $$
(with the convention $1/\infty:=0$).

\medskip
The extremal functions we shall encounter are obtained in the following
way. For  $j,k,l$ as above, 
it is well--known \cite[Vol.~I, p.~72]{Car} that there exists a
 hyperbolic triangle $\Delta$ in the unit disk $\D$ with 
interior angles $\pi/j$, $\pi/k$ and $\pi/l$. The triangle is moreover uniquely
determined up to a motion of the hyperbolic plane. 
The conformal map  from
$\Delta$ onto the upper halfplane $\H=\{w \in \C \, : \, \Im w>0\}$, which maps the vertex with angle $\pi/j$ to $0$,
the vertex with angle $\pi/k$ to $1$ and the vertex with angle $\pi/l$ to
$\infty$ is uniquely determined. By Schwarz reflection, this conformal map can be
analytically continued to a meromorphic function $f$ on $\D$ such that
all zeros of $f$ have exact order $j$, all zeros of $f-1$ have
exact order  $k$ and all poles of $f$ have  exact order $l$,
i.e., $f \in {\cal M}_{j,k,l}$. We call every such meromorphic
function a triangle map of order $(j,k,l)$. Note that a triangle map
of order $(j,k,l)$  is uniquely determined up to precomposition with a unit disk
automorphism. Clearly, a triangle map of order $(\infty,\infty,\infty)$ is a universal
covering from $\D$ onto $\C''$.

\begin{theorem}[Landau--type theorem] \label{thm:landau}
Let $j,k,l \ge 2$ be integers (or $=\infty$) such that $1/j+1/k+1/l<1$ and let
\begin{equation} \label{eq:const}
\begin{array}{rl} 
C_1& := \displaystyle j \arcsinh \left( \frac{1}{j \cdot \lambda_{1-1/j,1-1/k,1-1/l}(-1)} \right)
\, \\[4mm] 
\, C_3 & :=\displaystyle l \arcsinh \left( \frac{1}{l \cdot \lambda_{1-1/j,1-1/k,1-1/l}(-1)}
\right) \, . 
\end{array}
\end{equation}
Then for every $f \in  {\cal M}_{j,k,l}$ with  $a_0:=f(0)\not=\infty$, 
we have for $a_1:=f'(0)$ the sharp estimate
$$ |a_1| \le \begin{cases} 2 \,j \,|a_0| \, \sinh \left[ \displaystyle\frac{C_1+\big| \log |a_0|
      \big|}{j} \right] & \hspace{1cm}|a_0| \le 1\, \\
 & \text{ if } \\
2 \, l \, |a_0| \, \sinh \left[ \displaystyle \frac{C_3+\big| \log |a_0|
      \big|}{l}\right] & \hspace{1cm} |a_0| \ge 1\, .
\end{cases}$$
Equality holds if and only if $f$ is a triangle map
of order $(j,k,l)$ with $f(0)=-1$.
\end{theorem}


\begin{corollary}
Let $j,k \ge 2$ be integers (or $=\infty$) such that $2/j+1/k<1$. Then
for every $f \in {\cal M}_{j,k,j}$ with
 $a_0:=f(0)\not=\infty$, we have for $a_1:=f'(0)$ the sharp estimate
$$ |a_1| \le 2 \,j \,|a_0| \, \sinh \left[ \frac{C_1+\big| \log |a_0|
      \big|}{j} \right] \,  \, .$$
Here, $C_1$ is as in (\ref{eq:const}) with $l=j$, where  
$\lambda_{1-1/j,1-1/k,1-1/j}(-1)$ is given by (\ref{eq:f}) with
$\alpha_1=1-1/j$, $\alpha_2=1-1/k$ and $\alpha_3=1-1/j$.
 Equality holds if and only if $f$ is a triangle map
of order $(j,k,j)$ with $f(0)=-1$.
\end{corollary}

If $j=k=\infty$, then the corollary reduces to the well--known sharp 
version of Landau's theorem due to Hempel \cite{Hem79},
$$ |a_1| \le 2 |a_0| \Big( \big|\log |a_0|\big|+L \Big) \,  , \qquad
L=\frac{1}{\lambda_{1,1,1}(-1)}=\frac{1}{4 \pi^2} \cdot \Gamma \left( \frac{1}{4} \right)^4 \,  , $$
which holds for every analytic function $f(z)=a_0+a_1 z+\cdots$ in $\D$ 
 omitting $0$ and $1$. Equality occurs if and only if $f$ is a universal
covering from $\D$ onto $\C''$
 with $f(0)=-1$.

\begin{theorem}[Schottky--type theorem] \label{thm:schottky}
Let $j,k,l \ge 2$ be integers (or $=\infty$) such that $1/j+1/k+1/l<1$.
Then for every $f \in {\cal M}_{j,k,l}$ the sharp estimate
\begin{equation} \label{eq:s}
 \tanh \left( \frac{\tilde{C}_1+\log|f(z)|}{2 l} \right)
\le \tanh \left( \frac{\tilde{C}_1+\log^+|f(0)|}{2 l} \right)
\frac{1+|z|}{1-|z|} \, , \qquad z \in \D \, ,  
\end{equation}
holds, where $\log^+x=\max\{\log x,0\}$ for $x >0$ and
$$\tilde{C}_1=\displaystyle l \arcsinh \left( \frac{1}{l \cdot
    \lambda_{1-1/j,1-1/k,1-1/l}(-1)} \right)\, . $$ 
In particular,
 $$ \log |f(z)| \le 2\,  l \, \text{\rm arctanh} \left[ \tanh \left(
     \frac{\tilde{C}_1+\log^+|f(0)|}{2 l} \right) \frac{1+|z|}{1-|z|} \right]
 -\tilde{C}_1 $$
for all
$$ |z| < \exp \left(-  \frac{\tilde{C}_1+\log^+|f(0)|}{ l} \right) \, . $$
\end{theorem}

\begin{remark} \label{rem:sharp}
The estimate (\ref{eq:s}) is sharp in the following sense: if $M>0$
is a constant such that
$$ \tanh \left( \frac{M+\log|f(z)|}{2 l} \right)
\le \tanh \left( \frac{M+\log^+|f(0)|}{2 l} \right) \frac{1+|z|}{1-|z|} $$
holds for all $z \in \D$ and all
meromorphic functions $f \in {\cal M}_{j,k,l}$, then  one can show $M \ge \tilde{C}_1$.
\end{remark}


\begin{remark}
In the situation of Theorem \ref{thm:schottky}, we see that
if $f(0)\not=\infty$, then $$f(z)\not=\infty \text{ for all }
|z| <  \exp\left(-  \left(\tilde{C}_1+\log^+|f(0)|\right)/l \right) \, .$$
\end{remark}



If there are no poles $(l=\infty$) one gets a sharp Schottky--type result
on the entire unit disk:

\begin{corollary} \label{cor:schottky2}
Let $j,k \ge 2$ be integers (or $=\infty$) such that $1/j+1/k<1$. Then for
every $f \in {\cal M}_{j,k,\infty}$,
$$ \log |f(z)| \le \left[ C+\log^+|f(0)| \right] \frac{1+|z|}{1-|z|}-C \, ,
\qquad z\in\D \, , $$
where $ C=1/\lambda_{1-1/j,1-1/k,1}(-1)$.
\end{corollary}

\begin{corollary} \label{cor:schottky3}
Let $k \ge 2$ be an integer (or $=\infty$) and let 
$$ L_k:=\frac{1}{4 \pi^2} \cdot \Gamma \left( \frac{1+1/k}{4} \right)^2
\cdot \Gamma \left( \frac{1-1/k}{4} \right)^2 \cdot \cos \left( \frac{\pi}{2k}
\right) \, . $$ 
If $f$ is analytic and zero--free in $\D$ such that $f(z)-1$ has only zeros of order $\ge
k$, then
$$ \log |f(z)| \le \left[ L_k+\log^+|f(0)|\right] \frac{1+|z|}{1-|z|}-L_k \, ,
\qquad z \in \D \, .$$
\end{corollary}

A remark similar to Remark \ref{rem:sharp} applies to Corollary
\ref{cor:schottky2} as well as to Corollary \ref{cor:schottky3}.
Thus Corollary \ref{cor:schottky2} and Corollary \ref{cor:schottky3}
are in some sense best possible. The special case $k=\infty$ of Corollary
\ref{cor:schottky3} is the recent result of Li and Qi \cite{LQ2007}.

\section{Preliminaries} \label{sec:preliminaries}

We first recall a number of facts about conformal pseudo--metrics. Some of the
material is discussed in more detail in \cite{Hei62,Min1982}.


\smallskip

If $G$ is a domain in the complex plane $\C$, then 
we can identify a conformal pseudo--metric $\lambda(z) \, |dz|$ with its conformal
density,  that is the function $\lambda : G \to [0,+\infty)$, which 
represents the pseudo--metric $\lambda(z) \, |dz|$ in local coordinates when using the identity
map as a chart.
For instance, if $\lambda_{\alpha_1,\alpha_2,\alpha_3}(z)\, |dz|$ is the generalized
hyperbolic metric on $\P\backslash \{0,1,\infty\}$ of order
$(\alpha_1,\alpha_2,\alpha_3)$, 
then the associated generalized hyperbolic density
$\lambda_{\alpha_1,\alpha_2,\alpha_3}$ is a positive function on $\C''=\C
\backslash \{ 0,1\}$.

\smallskip

We call an upper semicontinuous pseudo--metric $\lambda(z)\, |dz|$ on $G \subset \C$
an SK--metric\footnote{Heins  \cite{Hei62} introduced the concept of
  SK--metrics and established a theory of such metrics. Note that he used the
  upper bound $-4$  in  his definition of SK--metrics instead of $-1$ as we do.}
if its (generalized) Gauss  curvature $\kappa_{\lambda}(z)$, defined by  
\begin{equation*}
\kappa_{\lambda}(z):=-  \frac{\liminf \limits _{r \to 0} \frac{4}{r^2} \left( \frac{1}{2\pi}
\int \limits_{0}^{2\, \pi} \log \lambda(z+re^{it})\, dt- \log\lambda(z)
\right)}{\lambda(z)^2}\, ,
\end{equation*}
is bounded above by $-1$ at every $z \in G$ with $\lambda(z)>0$.
Note, if $\lambda(z) \, |dz|$ is a {\it regular} conformal metric, i.e.,
$\lambda$ is twice continuously differentiable and strictly positive  on
$G$, then 
$\kappa_{\lambda}(z)=- \Delta \log \lambda(z)/\lambda(z)^2$, where
$\Delta$ denotes the usual Laplace operator.

\smallskip

The Fundamental Theorem about SK--metrics is Ahlfors' lemma \cite{Ahl38,Hei62}.
It  says that the hyperbolic metric $\lambda_{\D}(z) \, |dz|$ on the unit disk $\D$,
$$\lambda_{\D}(z) \, |dz|:=\frac{2 \, |dz|}{1-|z|^2} \, ,$$
is the maximal SK--metric on $\D$, i.e., $\mu(z) \le \lambda_{\D}(z)$ for
all $z \in \D$ and every SK--metric $\mu(z) \, |dz|$ on $\D$.
Actually, $\lambda_{\D}(z) \, |dz|$ is the {\it unique} maximal 
SK--metric on $\D$. This follows from the following result.

\begin{lemma}[Heins \cite{Hei62}] \label{lem:equality}
Let $\mu(z) \, |dz|$ be an SK--metric on a domain $G \subseteq \C$ and
$\lambda(z) \, |dz|$ a regular conformal metric on $G$ with constant curvature
$-1$ such that $\mu \le \lambda$. Then either $\mu<\lambda$ or $\mu \equiv
\lambda$.
\end{lemma}

By definition, the generalized hyperbolic metric
$\lambda_{\alpha_1,\alpha_2,\alpha_3}(z) \, |dz|$ is a regular
conformal metric on $\C''$ with constant curvature $-1$.
In general, conformal metrics with constant curvature play a
distinctive role. This comes in part from the well--known and easily verified fact that 
the Schwarzian $S_{\lambda}$ of  a regular  conformal  metric $\lambda(z) \, |dz|$ 
on a domain $G \subseteq \C$, \begin{equation} \label{eq:schwarz}
S_{\lambda}(z):=2 \left[ \frac{\partial^2 \log \lambda}{\partial z^2}(z)
 -\left( \frac{\partial \log \lambda}{\partial z}(z) \right)^2 \right] \, ,
\end{equation}
is a holomorphic function in $G$
if and only if  $\lambda(z) \, |dz|$ has constant curvature there.
The following classical fact tells us that locally every regular metric
with constant curvature $-1$ comes from the hyperbolic metric
$\lambda_{\D}(z) \, |dz|$ on the unit disk $\D$:

\begin{theorem}[Liouville \cite{Lio1853}] \label{thm:liouville}
Let $G \subseteq \C$ be a simply connected domain
and $\lambda(z) \, |dz|$  a regular conformal  metric on $G$
with constant curvature $-1$. Then the following are true.
\begin{itemize}
\item[(a)]
There exists a holomorphic function $\varphi : G \to \D$ such
that
\begin{equation} \label{eq:liouville}
 \lambda(z) =\frac{2\, |\varphi'(z)|}{1-|\varphi(z)|^2}  \, , \qquad z \in G
 \, . 
\end{equation}
The function $\varphi$ can be found among all solutions $\Psi$ to the Schwarzian differential
equation \begin{equation} \label{eq:schwarz0}
 \left( \frac{\psi''(z)}{\psi'(z)} \right)'-\frac{1}{2} \left(
   \frac{\psi''(z)}{\psi'(z)}
\right)^2=S_{\lambda}(z) \, . 
\end{equation}
\item[(b)]
An analytic function $g: G \to \D$ satisfies 
$$\lambda(z)=\frac{2\, |g'(z)|}{1-|g(z)|^2} \, , \qquad z \in G \,,$$
if and only if $g=T\circ \varphi$, where $T$ is an automorphism of $\D$.
\end{itemize}
\end{theorem}

Our derivation of the explicit formula for $\lambda_{\alpha_1,\alpha_2,\alpha_3}$
in Theorem \ref{thm}, which will be given in Section \ref{sec:proofs}, depends
in an essential way on  Liouville's theorem. 
Part (b) will also be used to show that the theorems of Landau and
Schottky--type stated in \S \ref{sec:apps} are best possible.

\section{Proofs} \label{sec:proofs}

\subsection{The explicit formula for the generalized hyperbolic metric}

The proof of Theorem \ref{thm} is based on the following lemmas.

\begin{lemma} \label{lem1}
Let $0<\alpha_1,\alpha_2 <1$ and $0 < \alpha_3 \le 1$ such that
$\alpha_1+\alpha_2+\alpha_3>2$  and define $\alpha, \beta, \gamma$ by (\ref{eq:hyp1}).
Then the following representation formulas are valid.
\begin{itemize}
\item[(a)] In the slit disk $\D^-=\D \backslash (-1,0]$ we have
$$ \lambda_{\alpha_1,\alpha_2,\alpha_3}(z)= \frac{2 \, |\varphi'(z)|}{1-|\varphi(z)|^2} \, , \qquad z \in \D^{-} $$
with
$$ \varphi(z)=c_0 \, \frac{z^{1-\gamma} F(\alpha-\gamma+1,\beta-\gamma+1,2-\gamma;z)}{
F(\alpha,\beta,\gamma;z)}$$
for some constant $c_0>0$.
\item[(b)] In the slit disk $K^+_1(1)=\{z \in \C \, : \, |z-1|<1\} \backslash [1,2)$ we have
$$ \lambda_{\alpha_1,\alpha_2,\alpha_3}(z)=  \frac{2 \, |g'(z)|}{1-|g(z)|^2}
\, , \qquad z \in K^+_1(1) \, , $$
with 
$$ g(z)=c_1 \, \frac{(1-z)^{\gamma-\alpha-\beta}
  F(\gamma-\beta,\gamma-\alpha,\gamma-\alpha-\beta+1;1-z)}{F(\alpha,\beta,\alpha+\beta-\gamma+1;1-z)}$$
for some constant $c_1>0$.
\end{itemize}
\end{lemma}

\begin{lemma} \label{lem2}
The constant $c_0$ in Lemma \ref{lem1} has the value
\begin{equation} \label{eq:c}
 \sqrt{
\frac{\Gamma (1-\alpha ) \, \Gamma (1-\beta ) \, \Gamma (\alpha+1 -\gamma )
  \,  \Gamma (\beta+1 -\gamma ) }{\Gamma (\alpha ) \, \Gamma
   (\beta )  \, \Gamma (\gamma -\alpha ) \, \Gamma (\gamma
   -\beta )}} \cdot \frac{\Gamma(\gamma)}{\Gamma(2-\gamma)}\, . 
\end{equation}
\end{lemma}

\begin{remark}
The proof of Lemma \ref{lem1} will show that
 the functions $\varphi$ and $g$ in Lemma \ref{lem1} can be analytically
continued along any path in $\C''$. Thus 
the representation formulas for $\lambda_{\alpha_1,\alpha_2,\alpha_3}(z)$ in Lemma \ref{lem1}
clearly hold for any $z \in \C''$.
\end{remark}

\begin{corollary} \label{cor}
Let $0<\alpha_1,\alpha_2 <1$ and $0 < \alpha_3 \le 1$ such that
$\alpha_1+\alpha_2+\alpha_3>2$
 and define $\alpha, \beta, \gamma$ by (\ref{eq:hyp1}).
Then 
$$ \lambda_{\alpha_1,\alpha_2,\alpha_3}(z)=\frac{2 \, c \, (1-\alpha_1)}{|z|^{\alpha_1} |1-z|^{\alpha_2} \left\{|F(\alpha,\beta,\gamma;z)|^2 -c^2 |1-z|^{2-2\alpha_1} |F(\alpha-\gamma+1,\beta-\gamma+1,2-\gamma;z)|^2 \right\}} \, $$
with $c$ given by (\ref{eq:c}).
\end{corollary}

In order to prove Lemma \ref{lem1}, we first need to recall some well--known
 facts. Let $\lambda(z) \, |dz|$ be a regular conformal metric on a 
domain $G \subseteq \C$ with constant curvature $-1$.
If  an isolated boundary point $z_j \in \C$ of $G$ is a singularity of
order $\alpha_j$ of $\lambda(z) \, |dz|$, then $S_{\lambda}$ has a pole
of order $2$ at $z_j$ and 
$$ S_{\lambda}(z)=\frac{(2-\alpha_j) \, \alpha_j}{2
  (z-z_j)^2}+\frac{c_j}{z-z_j}+O(1) \quad \text{ as } z \to z_j \, .$$
See, for instance, \cite{Nit57,KR2007,KR2008}.
Now let  $z_j=\infty$ be an isolated singularity of order $\alpha_j$ of
$\lambda(z)\, |dz|$. This means, by definition, that 
 $\mu(z) \, |dz|=\lambda(1/z) \, |dz|/|z|^2$ has an
isolated singularity of order $\alpha_j$ at $z=0$. Hence
$S_{\lambda}(z)=S_{\mu}(1/z)/z^4$ and we thus see that  
$$ \lim \limits_{z \to \infty} z^2 S_{\lambda}(z)=\lim \limits_{z \to 0}
z^2 S_{\mu}(z)=\frac{(2-\alpha_j)\, \alpha_j}{2} \, . $$


This observation leads to the following lemma.

\begin{lemma} \label{lem:0}
Let $z_1, \ldots, z_{n-1}$ and $z_n=\infty$ distinct points on $\P$ and let
$\lambda(z)\, |dz|$ be a regular conformal metric on $\C\backslash \{z_1, \ldots, z_{n-1}
\}$ with constant curvature $-1$ and singularities of order
$\alpha_j$ at $z_j$. Then
 $$ S_{\lambda}(z)=\sum \limits_{j=1}^{n-1} \left( \frac{(2-\alpha_j) \, \alpha_j}{2 \, (z-z_j)^2}+\frac{\beta_j}{z-z_j} \right) \, $$
and
$$ S_{\lambda}(z)=\frac{(2-\alpha_n) \, \alpha_n}{2 \,
  z^2}+\frac{\beta_n}{z^3}+O(1/z^4) \, , \qquad z \to \infty \, , $$
with complex numbers $\beta_1, \ldots, \beta_n$.
\end{lemma}
The numbers $\beta_1, \ldots, \beta_n $ are called the {\it accessory parameters} of $\lambda(z) \, |dz|$. In view of the asymptotic behavior of $S_{\lambda}(z)$ at $z=\infty$, the accessory parameters are related by
$$ \sum \limits_{j=1}^{n-1} \beta_j=0 \, , \quad
\sum \limits_{j=1}^{n-1} \left( (2-\alpha_j )\, \alpha_j+2 \beta_j z_j\right)=
(2-\alpha_n) \, \alpha_n \, , \quad \sum \limits_{j=1}^{n-1}
\left( (2-\alpha_j) \, \alpha_j z_j+\beta_j z_j^2 \right)=\beta_n \, . $$

In case of three singularities, these relations determine the accessory parameters completely.
Thus, if $\lambda_{\alpha_1, \alpha_2,\alpha_3}(z)$ is the generalized
hyperbolic density of order $(\alpha_1,\alpha_2,\alpha_3)$ on $\C''$, then
 $S_{\lambda_{\alpha_1,\alpha_2,\alpha_3}}$ is a rational
function with poles of order $2$ at $z=0$ and $z=1$ and
\begin{equation} \label{eq:schwarzhyp}
 S_{\lambda_{\alpha_1,\alpha_2, \alpha_3}}(z)=\frac{1}{2} \left[ \frac{1-\theta_1^2}{z^2}+\frac{1-\theta_2^2}{(1-z)^2}+\frac{1-\theta_1^2-\theta_2^2+\theta_3^2}{z \, (1-z)} \right]\, 
\end{equation}
with $ \theta_j=1-\alpha_j$, $j=1,2,3$. Hence in this case  the Schwarzian $S_{\lambda_{\alpha_1,\alpha_2, \alpha_3}}(z)$
is explicitly determined by $\alpha_1,\alpha_2,\alpha_3$.
In order to determine $\lambda_{\alpha_1,\alpha_2,\alpha_3}(z) $ from  $\alpha_1,\alpha_2,\alpha_3$,
we therefore need to recover the metric from its Schwarzian. Away from the
singularities one can use Theorem \ref{thm:liouville}
for this purpose. Thus we
 have to examine the Schwarzian differential equation
\begin{equation} 
\label{eq:schwarzhyp2}
 \left( \frac{\psi''(z)}{\psi'(z)} \right)'-\frac{1}{2} \left( \frac{\psi''(z)}{\psi'(z)}
\right)^2 =\frac{1}{2} \left[
  \frac{1-\theta_1^2}{z^2}+\frac{1-\theta_2^2}{(1-z)^2}+\frac{1-\theta_1^2-\theta_2^2+\theta_3^2}{z
    \, (1-z)} \right]\, . 
\end{equation}
and  use the following classical fact (see \cite[p.~116 ff.]{Car}).

\begin{lemma} \label{lem5}
Let $u_1,u_2$ be  two linearly independent solutions
of the hypergeometric differential equation
\begin{equation} \label{eq:hypergeometric}
 z (1-z) \, u''+\left[ \gamma-(\alpha+\beta+1) z \right] \, u'-\alpha \beta \,u=0 \, 
\end{equation}
with $\alpha$, $\beta$ and $\gamma$ determined by (\ref{eq:hyp1}).
Then the solutions $\psi$ of the  Schwarzian differential equation (\ref{eq:schwarzhyp2})
have the form $\psi(z)=T(u_2(z)/u_1(z))$ where $T$ is an arbitrary M\"obius transformation.
\end{lemma}

{\bf Proof of Lemma \ref{lem1}.}
We only prove part (a). The proof of part (b) is similar and is left to the 
reader. We consider (\ref{eq:hypergeometric}) in $\D^-$
and note that
$$ u^0_1(z)=F(\alpha,\beta,\gamma;z) \, ,  \qquad u^0_2(z)=z^{1-\gamma}
F(\alpha-\gamma+1,\beta-\gamma+1,2-\gamma;z)$$
are two linearly independent solutions to (\ref{eq:hypergeometric}).
 In view of Lemma \ref{lem5} and Theorem \ref{thm:liouville} (a), we know
that in $\D^-$ 
$$ \lambda(z):=\lambda_{\alpha_1,\alpha_2,\alpha_3}(z)=\frac{2 |\varphi'(z)|}{1-|\varphi(z)|^2} \quad \text{ with } \quad 
\varphi(z)=\frac{a u_2^0(z)+b u_1^0(z)}{c u_2^0(z)+d u_1^0(z)} $$
for appropriate constants $a,b,c,d \in \C$ with $ad-bc \not=0$.
If we
let $$h(z):=F(\alpha-\gamma+1,\beta-\gamma+1,2-\gamma;z)/F(\alpha,\beta,\gamma;z)
\, , $$
then a straightforward computation gives
\begin{equation} \label{eq:lambdi}
 |z|^{\gamma} \lambda(z)=\frac{2 \, |ad-bc| \, |(1-\gamma)h(z)+h'(z)|}{
\left[ |c|^2-|a|^2 \right]|z|^{2-2\gamma} |h(z)|^2+\left[ |d|^2-|b|^2
\right]+2 \Re \left[ \left(a\overline{b}-c\overline{d}\right) z^{1-\gamma} h(z) \right]}
\, .
\end{equation}
Since $h$ is analytic at $z=0$ with $h(0)=1$ and $|z|^{\gamma} \lambda(z)$ is
single--valued in $\D \backslash \{ 0 \}$, a glance at (\ref{eq:lambdi})
shows that $a \overline{b}=c \overline{d}$. Since $|z|^{\gamma} \lambda(z)$ is
strictly positive and $0<\gamma<1$, we can then deduce from (\ref{eq:lambdi})
that $|d| \ge |b|$. Moreover, we can exclude the case $|d|=|b|$, since
$\lambda(z) \, |dz|$ has a corner of order $\gamma$ at $z=0$, so
$|z|^{\gamma} \lambda(z)$ is bounded at $z=0$. Thus, $|b|<|d|$. Postcomposing
$\varphi$ with a unit disk automorphism $T$ which sends $b/d \in \D$ to $0$
and using Theorem \ref{thm:liouville} (b), we can hence assume that $b=0$ and
thus also $c=0$.
This proves part (a) with $c_0=a/d$. Note that we can take $c_0>0$ by
multiplying $\varphi$ with an appropriate complex number of absolute value
one. Thus $\varphi(z)=c_0\, u_2^0(z)/u_1^0(z)$ as claimed.
\hfill{$\blacksquare$}

\medskip

{\bf Proof of Lemma \ref{lem2}.}
Let again
$$ u^0_1(z)=F(\alpha,\beta,\gamma;z) \, ,  \qquad u^0_2(z)=z^{1-\gamma}
F(\alpha-\gamma+1,\beta-\gamma+1,2-\gamma;z)$$
be a fundamental system of (\ref{eq:hypergeometric}) in $\D$ and let
$$ u^1_1(z)=F(\alpha,\beta,\alpha+\beta-\gamma+1;1-z) \, ,  \qquad u^1_2(z)=(1-z)^{\gamma-\alpha-\beta}
F(\gamma-\beta,\gamma-\alpha,\gamma-\alpha-\beta+1;1-z)$$
be a fundamental system of (\ref{eq:hypergeometric}) in $K_1(1)$.
Note that the above fundamental systems are connected by the transition relations
\begin{equation} \label{eq:trans}
 u^0_1(z)=A u_1^1(z)+B u^1_2(z) \, , \qquad u^0_2(z)=C u_1^1(z)+D u_2^1(z) \, ,
\end{equation}
where 
$$ \begin{array}{rclcrll} \displaystyle 
A &=&  \displaystyle  \frac{\Gamma (\gamma ) \, \Gamma (\gamma-\alpha -\beta  )}{\Gamma (\gamma -\alpha )\,
   \Gamma (\gamma -\beta )} \, , \qquad &   
B &=&  \displaystyle  \frac{\Gamma(\gamma) \, \Gamma (\alpha +\beta -\gamma )}{\Gamma (\alpha )\, 
   \Gamma (\beta )} \, , \\[6mm]
C &=&  \displaystyle  \frac{\Gamma (2-\gamma ) \, \Gamma (\gamma-\alpha -\beta )}{\Gamma (1-\alpha )\,
   \Gamma (1-\beta )} \, , \qquad & 
D &=&  \displaystyle 
\frac{\Gamma (2-\gamma )\,  \Gamma (\alpha +\beta -\gamma )}{\Gamma (\alpha
   -\gamma +1) \, \Gamma (\beta -\gamma +1)} \, ,
\end{array}
$$
see \cite[Vol.~II, p.~141]{Car}.
By part (b) of Liouville's Theorem \ref{thm:liouville}
and Lemma~\ref{lem1}, we get
\begin{equation} \label{eq:1}
 \varphi(z)=\eta \, \frac{g(z)-z_0}{1-\overline{z_0} \, g(z)} \, .
\end{equation}
Here, $\varphi=c_0 u_2^0/u_1^0$ and $g=c_1 u_2^1/u_1^1$ are the functions of Lemma \ref{lem1}. Inserting these expression into (\ref{eq:1}) and using the transition relations (\ref{eq:trans}), we obtain
$$ c_1 e^{i t} \frac{\frac{u^1_2(z)}{u^1_1(z)}-\frac{z_0}{c_1}}{1-\overline{z_0} \, c_1 \frac{u^1_2(z)}{u^1_1(z)}}=\eta \, \frac{g(z)-z_0}{1-\overline{z_0} \, g(z)}=\varphi(z)=c_0 \frac{u_2^0(z)}{u^0_1(z)}=c_0 \frac{C u_1^1(z)+D u_2^1(z)}{A u_1^1(z)+B u^1_2(z)}=\frac{c_0D}{A} \frac{\frac{u^1_2(z)}{u_1^1(z)}+\frac{C}{D}}{1+\frac{B}{A} \frac{u_2^1(z)}{u_1^1(z)}} \, .
$$
This leads to
$$ c_1 e^{i t}= \frac{c\, D}{A} \, , \qquad -\frac{z_0}{c_1}=\frac{C}{D} \, , \qquad -\overline{z_0} \, c_1=\frac{B}{A} \,  $$
and therefore we get
$$ c_0=\sqrt{\frac{A B}{CD}} \,  .$$
An easy computation finally yields (\ref{eq:c}).
\hfill{$\blacksquare$}

\medskip

{\bf Proof of Corollary \ref{cor}.}
For
$$ \varphi(z)=c_0 \frac{z^{1-\gamma} F(\alpha-\gamma+1,\beta-\gamma+1,2-\gamma;z)}{
F(\alpha,\beta,\gamma;z)}$$
we have (\cite[Vol.~II, p.~147]{Car})
$$ \varphi'(z)= \frac{c_0 \, (1-\alpha_1)}{z^{\alpha_1} (1-z)^{\alpha_2} F(\alpha,\beta,\gamma;z)^2} \, ,$$
which proves the assertion of Corollary \ref{cor}. \hfill{$\blacksquare$}

\medskip

{\bf Proof of Theorem \ref{thm}}.
In order to prove  Theorem \ref{thm} we just use the representation formula of Corollary \ref{cor} and express
$u^0_2(z)=z^{1-\gamma}
F(\alpha-\gamma+1,\beta-\gamma+1,2-\gamma;z)$
in terms of
$u^0_1(z)=F(\alpha,\beta,\gamma;z)$  and $u_1^1(z)=F(\alpha,\beta,\alpha+\beta-\gamma+1;1-z)$
with the help of the transition formulas (\ref{eq:trans}).
\hfill{$\blacksquare$}

\medskip

We end this section with the following mapping properties of the function
$\varphi$ in Lemma \ref{lem1}~(a).

\begin{remark} \label{rem:f}
Let $\alpha_1, \alpha_2, \alpha_3$ be real parameters satisfying
condition (\ref{eq:gaussbonnet})
 and define $\alpha, \beta, \gamma$ by (\ref{eq:hyp1}).
Then the function
$$ \varphi(z)=c_0 \, \frac{z^{1-\gamma} F(\alpha-\gamma+1,\beta-\gamma+1,2-\gamma;z)}{
F(\alpha,\beta,\gamma;z)}$$
with $c_0$ given by (\ref{eq:c}) has an analytic continuation to  $\H$, which maps $\H$ conformally onto a hyperbolic triangle $\Delta
\subseteq \D$ with interior angles $\pi (1-\alpha_1)$, $\pi(1-\alpha_2)$ and $\pi (1-\alpha_3)$
 in such a way that $0$, $1$ and $\infty$ are mapped to the
vertices of $\Delta$. We refer to \cite[Vol.~II, p.~116 ff.]{Car}
for details.
\end{remark}

\subsection{Sharp lower bounds for the generalized hyperbolic density}

In order to prove the sharp lower bound (\ref{eq:lowerbound}) for the generalized
hyperbolic density $\lambda_{\alpha_1,\alpha_2,\alpha_3}$ we first state a simple, but important extremality property of
$\lambda_{\alpha_1,\alpha_2,\alpha_3}$.

\smallskip



\begin{lemma}\label{lem:unique2}
Let $\alpha_1, \alpha_2, \alpha_3$ be real parameters satisfying
(\ref{eq:gaussbonnet}).  
If $\mu(z) \, |dz|$ is an SK--metric on $\C''$ with singularities of order
$\beta_1 \le \alpha_1$, $ \beta_2 \le \alpha_2$ and $\beta_3 \le \alpha_3$ at
$z=0$, $z=1$ and $z=\infty$, then
$\mu \le \lambda_{\alpha_1,\alpha_2,\alpha_3}$.
Moreover, 
$\lambda_{\alpha_1,\alpha_2,\alpha_3}(z)\, |dz|$ is the unique conformal
metric on $\C''$ with constant curvature $-1$ and singularities of order
$\alpha_1$, $ \alpha_2$ and $\alpha_3$ at $z=0$, $z=1$ and $z=\infty$.
\end{lemma}

{\bf Proof.}
Let  $\lambda(z)\, |dz|$ be a conformal metric on $\C''$ with constant curvature
 $-1$ and singularities of order
$\alpha_1$, $ \alpha_2$ and $\alpha_3$ at $z=0$, $z=1$ and $z=\infty$. Then
the function 
$s(z)=\log^+\left(\mu(z)/\lambda(z)  \right)$
is subharmonic on $\C''$ in view of the curvature assumptions on $\mu$ and
$\lambda$. Moreover,  $s$ is bounded above at $z=0$ and $z=1$, so  it has a
subharmonic extension to $\C$. Since $u$ is also bounded above at $\infty$, we see
that $s\equiv c$ for some nonnegative constant. If $c>0$, then 
$\mu(z)=e^c\, \lambda(z)$, so $\kappa_{\mu}=e^{-2c} \kappa_{\lambda}>-1$, 
which violates the fact that $\mu(z) \, |dz|$ is an SK--metric.
 Hence $c=0$, so $\mu(z)\le\lambda(z)$ for all $z \in \C''$. 
Choosing $\lambda=\lambda_{\alpha_1,\alpha_2,\alpha_3}$ proves the first
part
of Lemma \ref{lem:unique2} and choosing
$\mu=\lambda_{\alpha_1,\alpha_2,\alpha_3}$ proves the second part.
 \hfill{$\blacksquare$}

\medskip

Hence 
$\lambda_{\alpha_1,\alpha_2,\alpha_3}(z)\, |dz|$ is  {\it maximal} 
among  all SK--metrics on
$\C''$ with singularities of order $\beta_1 \le \alpha_1, \beta_2\le \alpha_2$
and $\beta_3\le \alpha_3$ at $z=0,1$ and
$\infty$.
In order to make use of this maximality, we need the following simple 
``gluing lemma'' which we state for general SK--metrics.

\begin{lemma}[Gluing Lemma]\label{glueing}
Let $\lambda(z)\, |dz|$  be an SK--metric on  a domain $G\subset \C$
 and let $\mu(z)\, |dz|$ be an SK--metric on a subdomain $U$ of $G$ such that
the ``gluing condition''
\begin{equation*} 
\limsup \limits_{ U \ni z \to \xi} \mu(z) \le  \lambda(\xi)  
\end{equation*}
holds for all $\xi \in \partial U \cap G$. Then $\sigma(z)\, |dz|$ defined by
$$\sigma(z):=\begin{cases} \,   \max \{\lambda(z), \mu(z)\}   & \hspace{3mm} \, \text{for }  z \in  U \, , \\[2mm]
                       \,        \lambda(z)         & \hspace{3mm} \, \text{for } z \in G \backslash U
          \end{cases} $$
is an SK--metric on $G$.
\end{lemma}

{\bf Proof.}
The gluing condition guarantees that
 $\sigma$ is upper semicontinuous on $G$ and it is easy to see that
$\max \{\lambda(z), \mu(z)\}$ is the density of  an SK--metric on $U$. Hence
the curvature of $\sigma(z)\, |dz|$ is
 bounded above by $-1$ at
 each $z \in U$. If $z \in G\backslash U$ then $\kappa_{\sigma}(z)\le-1$.
This is clear if $z \not \in \partial U \cap G$. For $z \in \partial U \cap
G$, this follows from  $\sigma\ge \lambda$.
\hfill{$\blacksquare$}

\medskip

We now combine the maximality of $\lambda_{\alpha_1,\alpha_2,\alpha_3}(z) \,
|dz|$ with this gluing lemma.

\begin{theorem}[Strict Monotonicity]\label{thm:mono}
Let $\alpha_1,\alpha_2,\alpha_3$ be real parameters satisfying condition
(\ref{eq:gaussbonnet}).
Then $\lambda_{\alpha_1,\alpha_2,\alpha_3}(re^{it})$ is strictly decreasing for
$0<t<\pi$ and strictly increasing for $-\pi<t<0$ for any fixed $r\in(0,+\infty)$.
\end{theorem}

We note that the case $\alpha_1=\alpha_2=1$ of Theorem \ref{thm:mono}
was proved before by Hempel \cite{Hem79} if  $\alpha_3=1$ and by 
Anderson, Sugawa, Vamanamurthy and Vuorinen
\cite{A} if $\alpha_3<1$. The proofs in \cite{Hem79,A} are based on an a--priori
knowledge of the asymptotic behaviour of the metric at the corners, whereas
the following proof is solely based on the gluing lemma and the maximality of
the generalized hyperbolic metric.

\medskip

{\bf Proof of Theorem \ref{thm:mono}}. 
For  $\eta \in \partial \D$  let $\lambda_{\eta}(z)\, |dz|
:=\lambda_{\alpha_1,\alpha_2,\alpha_3}(\overline{\eta} \, z) \, |dz|$.
Then $\lambda_{\eta}(z)\, |dz|$ is the maximal SK--metric on 
 $\C \backslash \{ 0,\eta\}$
 with singularities of order $\alpha_1$, $\alpha_2$ and $\alpha_3$ at $z=0$, $z=\eta$ and $z= \infty$.
In a first step we show that
$\lambda_{\overline{\eta}}(z)=\lambda_{\eta}(\overline{z})$ for all $z\in
\C\backslash \{ 0 , \eta\}$. For this we note that 
$\lambda_{\eta}(\overline{z}) \, |dz|$ is an SK--metric 
 on $\C \backslash \{ 0 , \overline{\eta} \}$ with 
singularities of order $\alpha_1$, $\alpha_2$ and $\alpha_3$ at $z=0$,
$z=\overline{\eta}$ and $z= \infty$. Thus by maximality
\begin{equation}\label{eq:1a}
\lambda_{\eta}(\overline{z})
\le 
\lambda_{\overline{\eta}}(z) \, , \quad z \in \C \backslash \{0,  \overline{\eta} \}
\, .\end{equation}
Hence  $\lambda_{\overline{\eta}}(\overline{z}) \le
\lambda_{\eta}(z)$ for all $z \in \C \backslash \{0, \eta\}$ which implies 
\begin{equation}\label{eq:1b}
\lambda_{\overline{\eta}}(z)\le\lambda_{\eta}(\overline{z})\, ,\quad z \in \C
\backslash \{0,  \overline{\eta} \}\, .
\end{equation}
Combining (\ref{eq:1a}) and (\ref{eq:1b}) gives the desired result.

\smallskip

Second, we prove that $\lambda_{\overline{\eta}}(z) < \lambda_{\eta}(z)$ for all
  $z \in \H$ if $\Im \eta >0$. To check this assertion we note that 
$\lambda_{\overline{\eta}}(z)=\lambda_{\eta}(z)$ for all $z \in \R \backslash
\{0 \}$. Thus by the gluing lemma (Lemma \ref{glueing})
$$\sigma(z):= \begin{cases}
\max \{ \lambda_{\overline{\eta}}(z), \lambda_{\eta}(z) \} \, , &\quad z \in \H
\backslash \{\eta\}\\[2mm]
 \lambda_{\eta}(z) &\quad z \in \C \backslash ( \H \cup \{0 \}) 
\end{cases}
$$
induces an SK--metric on $\C \backslash \{0, \eta\}$ with
singularities of order $\alpha_1$, $\alpha_2$ and $\alpha_3$
at $z=0$, $z=\eta$ and $z=\infty$. Hence $\sigma \le \lambda_{\eta}$ and so
  $\lambda_{\overline{\eta}}(z) \le
\lambda_{\eta}(z)$ for all $z \in \H$ if $\Im \eta >0$ and Lemma
\ref{lem:equality} shows that  $\lambda_{\overline{\eta}}(z) <
\lambda_{\eta}(z)$ for all $z \in \H$ if $\Im \eta >0$.

\smallskip

Finally we derive the strict monotonicity of $\lambda(z)\,
|dz|:=\lambda_{\alpha_1,\alpha_2, \alpha_3}(z)\, |dz|$. Choose $\varphi_1,
\varphi_2 \in (- \pi , 0)$ with $\varphi_2 > \varphi_1$ and set $\eta_1:=
 e^{-i \, \varphi_1/2}$ and $\eta_2:= e^{i\,  \varphi_2/2}$. Then we have
\begin{equation*}
\begin{split}
\lambda(-r e^{i\varphi_1})&=\lambda(- r \, \overline{\eta_1}^2)=\lambda_{\eta_1}(-
r \overline{\eta_1})= \lambda_{\eta_1}(-r \overline{\eta_1} \, \eta_2 \,
\overline{\eta_2})= \lambda_{\eta_1\, \eta_2}(-r \overline{\eta_1} \, \eta_2)\\
&> \lambda_{\overline{\eta_1}\, \overline{\eta_2}}(-r \overline{\eta_1} \, \eta_2)
=\lambda(- r \eta_2^2)=\lambda(-r e^{i\varphi_2})%
\end{split}
\end{equation*}
~\hfill{$\blacksquare$}

\medskip

We are now prepared to prove Theorem \ref{thm:lowerbound}. We
shall use Theorem \ref{thm:mono} and one more time the gluing lemma.

\medskip

{\bf Proof of Theorem \ref{thm:lowerbound}.}
For $\alpha \le  1$ and $R>0$ let
$$\lambda_{\alpha,R}(z):=\frac{2 (1-\alpha) R^{1-\alpha}
  |z|^{-\alpha}}{R^{2(1-\alpha)}-|z|^{2(1-\alpha)}}=\frac{1-\alpha}{|z| \sinh
  \left[ (1-\alpha) \log \frac{R}{|z|} \right]} \, . $$
%
Here again, for the case $\alpha=1$ this formula has to be interpreted
in the limit sense $\alpha \nearrow 1$, i.e.,
$$ \lambda_{1,R}(z)=\lim \limits_{\alpha \nearrow 1} \lambda_{\alpha,R}(z)=
\frac{1}{|z| \log \frac{R}{|z|}} \, . 
$$
Then $\lambda_{\alpha,R}(z) \, |dz|$ is a conformal metric on the
punctured disk $0<|z|<R$ with constant curvature $-1$ and singularity
of order $\alpha$ at $z=0$. In point of fact, $\lambda_{\alpha,R}(z) \, |dz|$
is the maximal conformal metric on $0<|z|<R$ with those properties.

\smallskip
We now write $\lambda(z) \, |dz|:=
\lambda_{\alpha_1,\alpha_2, \alpha_3}(z) \, |dz|$. 
Note that $\lambda(z) \ge \lambda(-1)$ for all $|z|=1$ by Theorem
\ref{thm:mono}. If we choose $R_1$ such that $\lambda_{\alpha_1,R_1}(z)=
\lambda(-1)$ for $|z|=1$, i.e., $R_1:=e^{C_1}>1$, then
$$ \lambda(z)\ge \lambda(-1) =\frac{1-\alpha_1}{\sinh \left[ (1-\alpha_1)
    C_1\right]}=\frac{1-\alpha_1}{|z|\, \sinh \left[ (1-\alpha_1) \log \frac{R_1}{|z|}
\right]}=
\lambda_{\alpha_1,R_1}(z) \quad \text{ for all } \, |z|=1 \, . $$

So the gluing lemma (Lemma \ref{glueing}) ensures that
$$ \sigma(z):=\begin{cases}
\max\{ \lambda(z), \lambda_{\alpha_1,R_1}(z) \}\,  & \quad \text{if } \, 
0<|z|<1\, , \\
\lambda(z)\,  & \quad \text{if } \, |z| \ge 1 \, ,
\end{cases}$$
induces an SK--metric  on $\C''$ with  singularities of order
$\alpha_1$, $\alpha_2$ and $\alpha_3$ at $z=0$, $1$ and $\infty$.
The maximality of  $\lambda(z) \, |dz|$  implies
$\sigma(z) \le \lambda(z)$ for all $z \in \C''$. In particular,
$$\lambda(z) \ge \lambda_{\alpha_1,R_1}(z)=\frac{1-\alpha_1}{|z| \sinh \left[
    (1-\alpha_1) \log \frac{R_1}{|z|} \right]}=\frac{1-\alpha_1}{|z| \sinh \left[
    (1-\alpha_1) (C_1-\log|z|)  \right]}
$$ for all $|z| \le 1$, $z \not=0,1$,
with equality for $z=-1$. In a similar way,  one can prove
$$ \lambda(z) \ge  \frac{1-\alpha_3}{|z| \sinh \left[
    (1-\alpha_3) \log \left( R_3|z| \right)\right]}=\frac{1-\alpha_3}{|z| \sinh \left[
    (1-\alpha_3) (C_3+\log|z|)  \right]}$$
for all $|z| \ge 1$ with equality for $z=-1$. 

\medskip
Assume now there is $z_0 \in \C''$ such that equality holds in
(\ref{eq:lowerbound}). If $z_0 \in \D$, then
$\lambda(z_0)=\lambda_{\alpha_1,R_1}(z_0)$ and, as we have seen above,
$\lambda(z) \ge \lambda_{\alpha_1,R_1}(z)$ for all $0<|z|<1$.
Hence $\lambda(z) = \lambda_{\alpha_1,R_1}(z)$ for all $0<|z|<1$ by Lemma
\ref{lem:equality}. This however contradicts
$$\lim \limits_{z \to 1} \lambda_{\alpha_1,R_1}(z)<+ \infty=\lim \limits_{z \to 1}\lambda(z)\, .$$ 
In the same way, we can exclude the case $|z_0|>1$. Thus $|z_0|=1$, so $\lambda(z_0)=\lambda_{\alpha_1,R_1}(z_0)=\lambda(-1)$.
Now Theorem \ref{thm:mono} tells us that $z_0=-1$.
\hfill{$\blacksquare$}

\bigskip

{\bf Proof of Corollary \ref{cor:lowerbound}.}
Clearly, $C_1=C_3$ if $\alpha_1=\alpha_3$, so we only need to
compute the value $\lambda(-1)=\lambda_{\alpha_1,\alpha_2,\alpha_1}(-1)=
\lambda_{\alpha_1,\alpha_1,\alpha_2}(1/2)/4$.
Also,
Theorem \ref{thm} gives
$$ \lambda_{\alpha_1,\alpha_1,\alpha_2}(1/2)=
\frac{K_3}{1+K_1} \frac{2^{2 \alpha_1}}{\left| F\left(
        \alpha_1-\frac{\alpha_2}{2},\alpha_1-1+\frac{\alpha_2}{2},\alpha_1;\frac{1}{2}
      \right) \right|^2} \, . 
$$
Applying \cite[15.1.24]{Ab} and the duplication formula \cite[6.1.18]{Ab} for the
Gamma function, a straightforward computation gives
$$ F\left(
        \alpha_1-\frac{\alpha_2}{2},\alpha_1-1+\frac{\alpha_2}{2},\alpha_1;\frac{1}{2}
      \right)=\frac{2^{2 \alpha_1}}{8 \sqrt{\pi}} \frac{\Gamma(\alpha_1)
\Gamma \left( \frac{\alpha_1}{2}-\frac{\alpha_2}{4} \right) \Gamma \left(
  \frac{\alpha_1}{2}+\frac{\alpha_2}{4}-\frac{1}{2} \right)}{
\Gamma \left( \alpha_1-\frac{\alpha_2}{2} \right) \Gamma \left(
  \alpha_1+\frac{\alpha_2}{2}-1 \right)}\, . 
$$
On the other hand, using the reflection formula \cite[6.1.17]{Ab}
 in the expressions
(\ref{eq:consts}) for $K_1$ and $K_3$ and elementary trigonometric
manipulation lead to
$$ \frac{K_3}{1+K_1}=\sqrt{\frac{\tan \left( \frac{\pi}{2} \left(
        \frac{\alpha_2}{2}+\alpha_1 \right) \right)}{\tan \left( \frac{\pi}{2} \left(
        \frac{\alpha_2}{2}-\alpha_1 \right) \right)}} \frac{\Gamma(\alpha_1)^2}{\Gamma \left(
    \alpha_1-\frac{\alpha_2}{2} \right) \Gamma \left(
    \alpha_1-1+\frac{\alpha_2}{2} \right)} \, . $$
Combining the last two identities, we arrive at
$$ \lambda_{\alpha_1,\alpha_2,\alpha_1}(-1)=
\frac{\lambda_{\alpha_1,\alpha_1,\alpha_2}(1/2)}{4}=
\sqrt{\frac{\tan \left( \frac{\pi}{2} \left(
        \frac{\alpha_2}{2}+\alpha_1 \right) \right)}{\tan \left( \frac{\pi}{2} \left(
        \frac{\alpha_2}{2}-\alpha_1 \right) \right)}}
\frac{16 \pi}{2^{2 \alpha_1}} \frac{\Gamma \left(
    \alpha_1-\frac{\alpha_2}{2} \right)\Gamma \left(
  \alpha_1+\frac{\alpha_2}{2}-1 \right) }{\Gamma \left(
  \frac{\alpha_1}{2}-\frac{\alpha_2}{4} \right)^2 \Gamma \left(
  \frac{\alpha_1}{2}+\frac{\alpha_2}{4}-\frac{1}{2} \right)^2} \, .$$
Finally, making again use of the duplication formula \cite[6.1.18]{Ab} for the Gamma function 
for $z=\frac{\alpha_1}{2}-\frac{\alpha_2}{4}$ and
$z=\frac{\alpha_1}{2}+\frac{\alpha_2}{4}-\frac{1}{2}$ in the last numerator,
we deduce (\ref{eq:f}).
\hfill{$\blacksquare$}

\subsection{Schottky and Landau--type theorems}

We need the following variant of the Ahlfors--Schwarz lemma.

\begin{lemma}[Ahlfors--Schwarz]\label{lem:app1}
Let $j,k,l \ge 2$ be integers (or $=\infty$) such that $1/j+1/k+1/k <1$ and
let $\lambda(z) \, |dz|$ be the generalized hyperbolic density on $\C''$ 
of  order $(1-1/j,1-1/k,1-1/l)$. If $f \in {\cal M}_{j,k,l}$, then
$\lambda(f(z)) \, |f'(z)| \, |dz|$
is a regular conformal pseudo--metric of constant curvature $-1$ on $\D$, so 
$$ \lambda(f(z)) \, |f'(z)| \le \frac{2}{1-|z|^2}  \, , \qquad z\in \D\, . $$
 Equality for one point $z \in \D$ holds if and only
if $f$ is a triangle map of order $(j,k,l)$.
\end{lemma}

{\bf Proof.} Let $S:=f^{-1}(\{0,1,\infty\}) \cap \D$.
If $z_0 \in \D \backslash S$, then $\mu(z) \, |dz|:=\lambda(f(z)) \, |f'(z)| \,
|dz|$ is clearly a regular conformal pseudo--metric in a neighborhood of $z_0$
with constant curvature $-1$ there. Moreover, $\mu$ is 
continuous at any point $z_0 \in S$. In order to check this  for
the case $f(z_0)=0$, we note that the remainder function $r$ in
$$
\log \lambda(w) =  \begin{cases} -(1-1/j) \log |w|+r(w) & \text{ if } 2 \le j \le\infty  \\[3mm]
 - \log |w|-\log \left( -\log |w| \right)+r(w) & \text{ if } j=\infty \, .
\end{cases}$$
as $w \to 0$ is continuous at $w=0$.
This follows e.g.~from the results in \cite{KR2007}. Since
$f$ has a zero of order at least $j \ge 2$ at $z=z_0$, we easily deduce that
$\mu(z)=\lambda(f(z)) \, |f'(z)|$ is continuous at $z=z_0$.
The cases $f(z_0)=1$ and $f(z_0)=\infty$ are similar. Hence $\mu$ is
continuous on $\D$.
By \cite{Nit57}, $\mu(z) \, |dz|$ is actually regular on $\D \backslash \{z
\in \D \, : \, \mu(z)=0\}$ with curvature $-1$ there.
In particular, $\lambda(f(z)) \, |f'(z)|=\mu(z) \le \lambda_{\D}(z)$ for any $z \in
\D$ by the Ahlfors--Schwarz lemma. 
Note that if $f$ is a triangle map of order $(j,k,l)$, then by Remark
\ref{rem:f}, $\lambda(f(z)) \, |f'(z)|= \lambda_{\D}(z)$ for all $z \in \D$.

\smallskip
Conversely, if $\lambda(f(z)) \, |f'(z)|= \lambda_{\D}(z)$ 
 for some point $z \in
\D$, then $\lambda(f(z)) \, |f'(z)| \equiv \lambda_{\D}(z)$ in $\D$ by Lemma \ref{lem:equality}. Now pick a point $z_0 \in \D$ with $
w_0=f(z_0) \in f(\D) \backslash \{ 0,1,
\infty\}$. Then $f'(z_0)\not=0$, so $f$ has a local inverse $h$ in some disk
$K_r(w_0)$ such that $h(K_r(w_0)) \subset \D$. Hence
$\lambda(w)=\lambda_{\D}(h(w)) \, |h'(w)|$ for all $w\in K_r(w_0)$. 
Shrinking $r>0$ if necessary, we also have
$ \lambda(w)=\lambda_{\D}(\varphi_0(w)) \, |\varphi_0'(w)|$ in $K_r(w_0)$
where $\varphi_0$ is a local inverse of a triangle map $f_0$ of order $(j,k,l)$.
By Theorem \ref{thm:liouville} (b), we get
$\varphi_0=T \circ h$ for some disk automorphism $T$, so $f=f_0 \circ T$, i.e.,
$f$ is a triangle map of order $(j,k,l)$. 
\hfill{$\blacksquare$}

\medskip

{\bf Proof of Theorem \ref{thm:landau}.}
Let $\lambda(w) \, |dw|$ be the generalized hyperbolic density on $\C''$ with
singularities 
of order $(1-1/j,1-1/k,1-1/l)$.
Then $ \lambda(f(z)) \, |f'(z)| \le \lambda_{\D}(z)$ for each $z \in \D$ by
Lemma \ref{lem:app1}.
For $z=0$, we get 
 $|a_1|=|f'(0)| \le 2/\lambda(f(0))=2/\lambda(a_0)$.
Now employing the lower bound for $\lambda$ provided by Theorem
\ref{thm:lowerbound} gives the estimate of Theorem \ref{thm:landau}. 

\medskip

To handle the case of equality, we note 
Lemma \ref{lem:app1} and Lemma \ref{lem:equality} show that $f: \D \to \C$  is a triangle
map of order $(j,k,l)$  if and only if
$\lambda(a_0)\, |a_1|=2$. By Theorem \ref{thm:lowerbound} we have
$$\lambda(a_0)= \frac{1}{j \, |a_0|}\,
\frac{1}{\sinh\left[ \frac{1}{j}\, \left( C_1 + |\log|a_0||\right)   \right]}=
\frac{1}{l \, |a_0|}\,
\frac{1}{\sinh\left[ \frac{1}{l}\, \left( C_3 + |\log|a_0||\right)
  \right]}$$
if and only if $a_0=-1$. This finishes the proof.\hfill{$\blacksquare$}




\medskip

{\bf Proof of Theorem \ref{thm:schottky}.} 
 Let $g:=1/f \in {\cal M}_{l,k,j}$ and 
let $\lambda(z)$ denote the generalized hyperbolic
density on $\C''$ of order $(1-1/l,1-1/k,1-1/j)$. 
Then Lemma \ref{lem:app1} gives
\begin{equation} \label{eq:ahl}
\lambda(g(z)) \, |g'(z)| \le \frac{2}{1-|z|^2} \quad  \text{ for all } z \in \D \, ,
\end{equation}

Pick a point $z_0 \in \D$ such that
$|g(z_0)|<1$ and consider the curve $\gamma(t):=g(t \eta)$ for $t \in
[0,|z_0|]$ and $z_0=|z_0| \eta$. If $\gamma \subset \D$, then
Theorem \ref{thm:lowerbound} applied for $z=t \eta$ and (\ref{eq:ahl})
 lead to
\begin{equation} \label{eq:hihi}
 \frac{|g'(t\eta)|}{l \, |g(t\eta)| \, \sinh \left[ (\tilde{C}_1-\log |g(t \eta)|)/l \right]} \le
\frac{2}{1-t^2} \, , \qquad t \in [0,|z_0|] \, . 
\end{equation}
Integrating over
$[0,|z_0|]$ using $\frac{d}{dt} |g(t \eta)| \le |g'(t \eta)|$ yields
\begin{eqnarray*}
\int \limits_{|g(0)|}^{|g(z_0)|} \frac{|ds|}{l \, s \, \sinh\left[ (\tilde{C}_1-\log
    s)/l \right]} \le \log \frac{1+|z_0|}{1-|z_0|} \, . 
\end{eqnarray*}
Hence
\begin{equation} \label{eq:q1}
\left| \log \left[\frac{\tanh \left( \frac{\tilde{C}_1-\log |g(z_0)|}{2\,  l} \right)}{\tanh
    \left( \frac{\tilde{C}_1-\log |g(0)|}{2\, l} \right)} \right]\right|
 \le \log  \frac{1+|z_0|}{1-|z_0|} \, . 
\end{equation}
If $\gamma \not\subset \D$, then a similar argument using the ``last'' point
$\gamma(t^*)$ of
$\gamma$ outside $\D$ and integrating (\ref{eq:hihi}) from $t^*$ to $|z_0|$
gives
\begin{equation} \label{eq:q2}
 \tanh \left( \frac{\tilde{C}_1-\log |g(z_0)|}{2 \, l} \right) \le \left[\tanh\left(
  \frac{\tilde{C}_1}{2 \, l} \right) \right]   \cdot \frac{1+|z_0|}{1-|z_0|}\, . 
\end{equation}
Thus in both cases, $\gamma \subset \D$ and
$\gamma \not\subset \D$, we get by the monotonicity of $\tanh$
\begin{equation}\label{eq:15}
\tanh \left( \frac{\tilde{C}_1-\log |g(z_0)|}{2 l} \right) \le \left[\tanh\left(
  \frac{\tilde{C}_1+\log^+\frac{1}{|g(0)|}}{2 l} \right) \right]   \cdot
\frac{1+|z_0|}{1-|z_0|}\, . 
\end{equation}
If $|g(z_0)| \ge 1$, then (\ref{eq:15}) is trivially true.
Finally going back to $f=1/g$ finishes the proof.~\hfill{$\blacksquare$}

\medskip






\medskip

{\bf Proof of Corollary \ref{cor:schottky3}.}
Choosing $j=l=\infty$ in Theorem
 \ref{thm:schottky} gives $L_k=1/\lambda_{1,1-1/k,1}(-1)$. Equation (\ref{eq:f})
   shows
$$\lambda_{1,1-1/k,1}(-1)=2\, \, \frac{\Gamma \big(3/4+1/(4k)\big) \,\Gamma
    \big(3/4- 1/(4k)\big)}{\Gamma \big(1/4+1/(4k)\big) \,\Gamma
    \big(1/4- 1/(4k)\big)}\, .$$
Applying \cite[6.1.18]{Ab} for $z=1/4+1/(4k)$ and $z=1/4-1/(4k)$ in the
numerator and then using \cite[6.1.17]{Ab} gives the desired result.
\hfill{$\blacksquare$}


\end{document}